\documentclass[12pt,twoside,leqno]{article}
\usepackage{amsmath}
\usepackage{amssymb}
\usepackage{amsxtra}
\usepackage{amscd}
\usepackage{amsthm}
\usepackage[mathscr]{eucal}
\usepackage{amsmath,amssymb,latexsym,amsfonts,graphicx}
\usepackage{array}
\usepackage[all]{xy}
\usepackage{supertabular}

\usepackage[color,notcite,notref]{showkeys}
\definecolor{labelkey}{rgb}{1,0,0}

\newtheorem{thm}[subsection]{Theorem}
\newtheorem{prop}[subsection]{Proposition}
\newtheorem{cor}[subsection]{Corollary}
\newtheorem{lem}[subsection]{Lemma}

\theoremstyle{definition}

\input cyracc.def 
\font\tencyr=wncyr10
\def\sha{\text{\tencyr\cyracc{Sh}}}

\newcommand{\Q}{\mathbb Q}

\newcommand{\C}{\mathbb C}
\newcommand{\Z}{\mathbb Z}

\newcommand{\ddz}{\frac{d}{dz}}
\newcommand\hb{\hfil\break}

\newcommand\End{\mathrm{End}\,}
\newcommand\Ker{\mathrm{Ker}\,}

\newcommand\Gal{\mathrm{Gal}\,}

\newcommand\ord{\mathrm{ord}\,}

\begin{document}

\title{The Tate-Shafarevich group for elliptic curves with complex multiplication II}

\author{J. Coates, Z. Liang, R. Sujatha}


\maketitle

\section{Introduction}

Let $E$ be an elliptic curve over $\Q$ . Put $g_{E/\Q}=$ rank of $E/\Q,$ and
$$
\sha(E/\Q)=\Ker\left(H^1(\Q,E) \to \underset{v}{\oplus}\,H^1(\Q_v,E)\right),
$$
where $v$ runs over all places of $\Q$, with $\Q_v$ the completion of $\Q$ at $v$.
Although no algorithm has ever been  proven to work infallibly, the group $E(\Q)$ is,
in fact, easy to
determine in practice. By contrast, $\sha(E/\Q)$ is extremely difficult to study
either theoretically or numerically. The aim of the present paper
is to strengthen the theoretical and numerical results of \cite{CZS},
assuming $E$ has complex multiplication.

\medskip

For each prime $p$, let $t_{E/\Q,p}$ denote the $\Z_p$-corank of the $p$-primary subgroup of $\sha(E/\Q).$

\begin{thm}\label{1.1}
Assume that $E/\Q$ admits complex multiplication. For each $\epsilon>0$, there exists an explicitly computable number $c(E,\epsilon)$ such that
\begin{equation}\label{teq}
t_{E/\Q,p} \leq (1/2 +\epsilon)p - g_{E/\Q}
\end{equation}
for all $p \geq c(E,\epsilon)$  where $E$ has good ordinary reduction.
\end{thm}
\noindent  Of course, this result is a far cry from the standard conjecture that $t_{E/\Q,p}=0$ for every prime $p$. Our method of proof is similar to that given in \cite{CZS},  but we obtain the stronger result by employing an interesting observation of Katz \cite{K}  about the divisibility of the relevant $L$-values by primes where $E$ has good supersingular reduction. However, we stress that  Katz's work  is used to obtain information about the Iwasawa theory at good ordinary primes $p$, and we have no idea at present how to prove a result like \eqref{teq} for all sufficiently large primes $p$  where $E$ has good supersingular reduction.

Secondly, we extend the numerical computations of $t_{E/\Q,p}$ given in \cite{CZS} for certain $E$ with $g_{E/\Q}\geq 2$. Let
\begin{equation}\label{eqe1}
E\,:\, y^2=x^3-82x.
\end{equation}
Then  $g_{E/\Q}=3$ and $E(\Q)$ modulo torsion is generated by the points
$$
(-9,3), \,(-8,12),\, (-1,9).
$$
The conjecture of Birch and Swinnerton-Dyer predicts that $\sha(E/\Q)=0$ for this curve, but of course this is unproven, and we do not know whether $\sha(E/\Q)$ is even finite.

\begin{thm}\label{1.2}
For the elliptic curve \eqref{eqe1}, we have $\sha(E/\Q)(p)=0$ for all primes $p\equiv 1 \mod\! 4$ with $p \neq 41$ and $p< 30,000.$
\end{thm}
\noindent In fact,  a different and more subtle technique than that
used in \cite{CZS} is required to carry out these computations,
because this earlier  method  relies on calculating traces from the
field of 328-division points on the curve \eqref{eqe1}. This field
has the enormous degree 12,800 over $\Q(i)$ (the integer 328 occurs
here because the conductor of the Gr\"ossencharacter of the curve is
$328\Z[i]$), and we have been unable to find the minimal equation
over $\Q(i)$ of the $x$-coordinate of a 328-division point. By some
curious arguments in Galois theory (see Lemma \ref{ro1} and
\eqref{72}), we show that it suffices to work with a subfield of
degree 6,400 over $\Q(i)$, and,  with the help of MAGMA, we have
succeeded in explicitly computing the minimal polynomial of a
natural generator for this subfield over $\Q(i)$  (the
$x$-coordinate of a $328/(1+i)$-division point).

\medskip

Moreover, using the same technique of calculation, we have also carried out computations on the five curves of rank 2 given by
$$
E_i\,:\, y^2=x^3-D_ix,~~~{\rm with}~ D_1=-14,\, D_2=17,\,D_3=-33,\,D_4=-34,\, D_5=-39,
$$
extending the numerical results given in \cite{CZS} for the curves $E_1$ and $E_2.$

\medskip

\begin{thm}\label{ei}
For each of the curves $E_i ~(i=1,\cdots, 5),~ \sha(E_i/\Q)(p) $ is finite for all primes
$p$ of good reduction with $p \equiv 1 \mod\!4$ and $p < 30,000.$ Moreover, $\sha(E_i/\Q)(p)=0$
for all such $p$, except possibly in the four exceptional cases given by
$p=29$ and $277$ for the curve $E_1$,
$p=577$ for the curve $E_4$, and $p=17$ for the curve $E_5.$
\end{thm}

\noindent In fact the conjecture of Birch and Swinnerton-Dyer predicts that $\sha(E_i/\Q)=0$ for $i=1,\cdots,5.$  We ourselves have not successfully carried out the calculations of $p$-adic heights to verify that $\sha(E_i/\Q)(p)=0$ in the four exceptional cases of Theorem \ref{ei}  (our claim to have done this for the primes $p=29$ and $277$ for $E_1$ in \cite{CZS} is not correct). However, we are very grateful to C. Wuthrich for computing the p-adic heights for the three exceptional primes p=17, 29, and 277, thereby
confirming that  $\sha(E_i/\Q)(p)=0$ for the relevant curves for these primes  (the remaining exxceptional prime of 577 for the curve $E_4$ remains unsettled).

\medskip

Finally,  we would like to thank Allan Steel and Mark Watkins for their help in factoring  polynomials.

\section{Divisibility of $L$-values by supersingular primes}

As Katz has pointed out in \cite{K}, the special values of complex $L$-functions attached to elliptic curves with complex multiplication tend to be highly divisible by supersingular primes. In this section, we use his method to establish the precise version of his results  needed to prove Theorem \ref{1.1}.

\medskip

The following notation will be used throughout the rest of this paper . Let $K$ be an imaginary quadratic field embedded in the field $\C$ of complex numbers, and   ${\cal O}_K$   its ring of integers. Let $E$ be an elliptic curve defined over $K$ such that
$$
\End_K(E) \simeq {\cal O}_K.
$$
The existence of such a  curve implies that $K$ necessarily has class number 1. We fix a global minimal Weierstrass equation for $E$
\begin{equation}\label{eqe}
y^2+a_1xy+a_3y =x^3 +a_2x^2 +a_4 x+a_6,
\end{equation}
whose coefficients belong to ${\cal O}_K$.  Let $\psi_E$ denote the Gr\"ossencharacter of $E$ over $K$, and write ${\frak f}$ for the conductor of $\psi_E$. For each integer $n \geq 1$, define
$$
L_{\frak f}(\bar{\psi}_E^n,s)=\,\underset{(v,\frak f)=1}{\prod} \,\left( 1-\frac{\bar{\psi}_E^n(v)}{(Nv)^s} \right)^{-1}.
$$
This $L$-function is not, in general, primitive, and we write $L(\bar{\psi}_E^n,s)$ for the primitive Hecke $L$-function of $\bar{\psi}_E^n.$

\medskip

Let $\cal L$ be the period lattice of the N\'eron differential of the equation \eqref{eqe}, and let
\begin{equation}\label{pz}
\Phi(z,{\cal L})\,:\, \C/{\cal L} \simeq E(\C)
\end{equation}
be the isomorphism given by
$$
\Phi(z,{\cal L})=\left(\wp(z,{\cal L})-\frac{a_1^2+4a_2}{12},\, 1/2\left(\wp'(z,{\cal L})-
a_1\left(\wp(z,{\cal L})-\frac{a_1^2+4a_2}{12}\right)-a_3\right)\right),
$$
where $\wp(z,{\cal L})$ denotes the Weierstrass $\wp$-function attached to ${\cal L}.$ We also fix
$\Omega_{\infty}$ in ${\cal L}$ such that ${\cal L}=\Omega_{\infty}{\cal O}_K.$

\medskip

As we shall see in the proof of the next theorem, the numbers
\begin{equation}
\Omega_{\infty}^{-n}L(\bar{\psi}_E^n,n)    \, \, (n=1,2,...)
\end{equation}
all belong to $K$. Our goal is to use Katz's method to establish the following specific result. If $b$ is a real number, $[b]$
will denote as usual the largest integer $\leq b$.

\begin{thm}\label{kat}
Let $q$ be an odd prime number which is inert in $K$. Assume that $E$ has good reduction at $v=
q{\cal O}_K$. Then, for all integers $n\geq 3$, which are not congruent to $1+q$ modulo $(q^2-1)$, we have
\begin{equation}\label{super}
ord_v((n-1)!\Omega_{\infty}^{-n}L_{\frak f}(\bar{\psi}_E^n,n) )  \geq  { \left[\frac{nq}{(q^2-1)}\right]-1}.
\end{equation}
\end{thm}

\medskip

The proof of this theorem will take up the rest of this section. The initial arguments, although subsequently used for supersingular primes,  are  motivated by the well known construction of the $p$-adic
$L$-functions of $E$ in the ordinary case (see \cite{CW}, \cite{CS}). If $\alpha$  is any non-zero element of $ {\cal O}_K$,
let $E_{\alpha}$ denote the kernel of multiplication by $\alpha$ on $E(\bar{K})$. Fix for the rest of the paper an element $f$
of $ {\cal O}_K$ such that $\frak f = f {\cal O}_K$. The field
\begin{equation}\label{div}
F = K(E_f)
\end{equation}
will play an important role in  both our theoretical and numerical arguments. It is an abelian extension of $K$, which is ramified precisely at the bad primes of $E$ over $K$, and it coincides with the ray class field of $K$ modulo $\frak f$ (see, for example,
\cite{CS}, Chap. 2). Moreover, if $w$ is any good prime of $E/K$, we have
\begin{equation}\label{gr1}
{\psi_E}(w)(U) = U^{\tau_w}
\end{equation}
for all $U$ in $E_f$, where $\tau_w$ denotes the Artin symbol of $w$ for $F/K$. Put
\begin{equation}\label{gr2}
G = Gal(F/K).
\end{equation}

\medskip

Let $\lambda$  be any element of  ${\cal O}_K$, which is not a unit, and which is relatively prime to $6fq$. Let $J_{\lambda}$ denote the set of all non-zero elements in $E_f$. Define
\begin{equation}\label{r1}
R_{\lambda}(P) =   c_{E}(\lambda).\underset{U}{\prod}(x(P) - x(U))^{-1},
\end{equation}
where $ U$  runs over any  set of representatives of $J_{\lambda}$  modulo the action of the group $\mu_2$ of square roots of unity,  and $c_E(\lambda)$ is the unique non-zero element of $K$, whose existence is established in  \cite{LMS} (see Proposition 1 of the Appendix). Thus $R_{\lambda}(P)$ is a rational function on $E$ with coefficients in $K$. Let $V=\Phi(\Omega_{\infty}/f, \cal L)$,
and define
\begin{equation}\label{r2}
\frak R_{\lambda}(P) = \underset{\tau \in G}{\prod}R_{\lambda}(P \oplus V^{\tau}),
\end{equation}
where $\oplus$ denotes the group law on $E$. Thus $\frak R_{\lambda}(P)$ is also a rational function on $E$, with coefficients in $K$. In view of \eqref{gr1} and Theorem 3 of the Appendix of \cite{LMS}, we have
\begin{equation}\label{r3}
\frak R_{\lambda}({\psi_E}(w)(P) ) = \underset{S \in E_w}{\prod}\frak R_{\lambda}(P \oplus S)
\end{equation}
for each prime $w$ of $K$, where $E$ has good reduction for $E/K$. Hence, defining
\begin{equation}\label{psil}
\Psi_{\lambda}(P)={\frak R}_{\lambda}(P)^{Nv}/
{\frak R}_{\lambda}(\psi_E(v)(P)),
\end{equation}
it follows from \eqref{r3} with $w=v$ that
\begin{equation}\label{psilp}
\underset{S \in E_q}{\prod}\,\Psi_{\lambda}(P \oplus S)=1.
\end{equation}

\medskip

The following result shows that the function $\Psi_{\lambda}(P)$ is related to $L$-values. For each integer $n \geq 1$, put
\begin{equation}\label{ln}
L_n=(n-1)!\,\Omega_{\infty}^{-n}L_{\frak f}(\bar\psi_E^n,n).
\end{equation}

\begin{prop}\label{ddz}
For all integers $n \geq 1$, we have
\begin{equation}\label{dndz}
(Nv)^{-1}\left(\ddz\right)^n \log\, \Psi_{\lambda}(\Phi(z,{\cal L}))_{\mid_{z=0}}=(-1)^{n-1}f^n(N\lambda-\psi_E((\lambda))^n)\left(1-\frac{\psi_E(v)^n}{Nv}\right)L_n.
\end{equation}
\end{prop}
\proof
This is a very classical calculation (see for example \cite{CW}), and we only sketch the main points in the proof. Put
$\pi=\psi_E(v).$ Since $\pi(\Phi(z,{\cal L}))=\Phi(\pi z, {\cal L}),$ we conclude from \eqref{psil} that it suffices to prove that, for all integers $n \geq 1$,
\begin{equation}\label{ddzl}
\left(\ddz\right)^n\,\log {\frak R}_{\lambda}\left(\Phi(z,{\cal L})\right)_{\mid{z=0}}\,=\,(-1)^{n-1}f^n\left(N\lambda- \psi_E((\lambda))^n\right)L_n.
\end{equation}
Let
$$
\theta(z,{\cal L})=\exp\left(-s_2({\cal L})z^2/2\right)\sigma(z,{\cal L}),
$$
where $\sigma(z,{\cal L})$ is the Weierstrass $\sigma$-function of ${\cal L}$ and, as usual,
$$
s_2({\cal L})=\underset{\underset{s>0}{s \rightarrow 0}}{\lim}\,
\underset{w \in {\cal L}\setminus{0}}{\Sigma}\,
w^{-2}\mid w\mid^{-2s}.
$$
Then we have (see, for example, \cite[Theorem 1.9]{GS})
\begin{equation}\label{rlt}
\frak R_{\lambda}(\Phi(z,{\cal L}))^2=c_E(\lambda)^{2d} \,\underset{{\frak b} \in {\cal B}}{\prod}\, \theta^2(z+\psi_E({\frak b})
\frac{\Omega_{\infty}}{f}, {\cal L})^{N\lambda}/\theta^2(z+
\psi_E({\frak b})\frac{\Omega_{\infty}}{f}, \, \lambda^{-1}{\cal L}),
\end{equation}
where ${\cal B}$ is any set of integral ideals of $K$ whose Artin symbols for $F/K$ give precisely the elements of the Galois group $G$, and $d= [F:K]$.

\medskip

For each integer $n \geq 1$, let $E_n^*(z,{\cal L})$ be the value at $s=n$ of the analytic continuation of the Kronecker series
$$
H_n(z,s,{\cal L})\,=\,\underset{w \in \cal L}{\Sigma}\,\frac
{(\bar z+\bar w)^n}{\mid z+w \mid^{2s}}.
$$
As usual, let $A({\cal L})=(u\bar v-v\bar u)/2\pi i,$ where $u,v$ is any $\Z$-basis of ${\cal L}$ with $v/u$ having positive imaginary part. Then we have (see \cite[Corollary 1.7]{GS}), for any $\rho$ in $\C\setminus {\cal L},$
\begin{equation}\label{drh}
\ddz\, \log(\theta(z+\rho, {\cal L})) = \frac{\bar \rho}{A(\cal L)} + \overset{\infty}{\underset{n=1}{\Sigma}}\,(-1)^{n-1}
E_n^*(\rho,{\cal L}) z^{n-1}.
\end{equation}
Since
$$
A(\lambda^{-1}{\cal L})=A({\cal L})/N\lambda,~E^*_n(z,\lambda^{-1}{\cal L})\,=
\psi_E((\lambda))^n\,E_n^*\left(\psi_E((\lambda))z,{\cal L}\right),
$$
it follows easily from \eqref{rlt} and \eqref{drh}, on putting
$$ {\cal D}_{\lambda}=\left(\ddz\right)^n\log\,{\frak R}_{\lambda}
\left(\Phi(z,{\cal L})\right)_{\mid_{z=0}},
$$
that
\begin{equation}\label{drdz}
{\cal D}_{\lambda}
\,=\,(-1)^{n-1}\,(n-1)!\,
\underset{{\frak b} \in {\cal B}}{\sum}\,
\left(N\lambda\, E_n^*(\psi_E({\frak b})\frac{\Omega_{\infty}}{f},{\cal L})
- \psi_E((\lambda))^n\, E_n^*(\psi_E(\frak{b}\cdot(\lambda))
\frac{\Omega_{\infty}}{f},{\cal L})\right).
\end{equation}
On the other hand, it is easily seen  that, for all integers
$n \geq 1$, we have
\begin{equation}\label{lfp}
L_{f}(\bar \psi_E^n,s)\,=\,\frac{\mid\Omega_{\infty}/f\mid^{2s}}
{(\overline{\Omega_{\infty}/f)^n}}\,\underset{{\frak b}\in
{\cal B}}{\sum}\,H_n\left(\frac{\psi_E({\frak b})
\Omega_{\infty}}{f},s,{\cal L} \right).
\end{equation}
This completes the proof of \eqref{ddzl}, and so also of
Proposition {\ref{ddz}}
\qed

\medskip

We now turn to the $v$-adic properties of our function
$\Psi_{\lambda}(P).$ Let $\hat{E}^v$ be the formal group of $E$ at $v$, which is defined over ${\cal O}_v,$ and has parameter $t=-x/y.$ It can be shown that the action of ${\cal O}_K$ on $E$ extends to an action of ${\cal O}_v$ on
$\hat{E}^v.$ If $a$ is any element of ${\cal O}_v$, we write $[a](t)$ for the formal power series in ${\cal O}_v[[t]]$
giving the corresponding endomorphism of $\hat{E}^v.$ As before, put $\pi=\psi_E(v)$. Since the reduction modulo $v$ of the endomorphism $\pi$ of $E$ gives the Frobenius endomorphism of the reduced curve, we plainly have
\begin{equation}\label{pit}
[\pi](t) \equiv\, t^{Nv}\!\mod q.
\end{equation}
As $[a](t)=at+\cdots$ for all $a \in {\cal O}_v,$ it follows
from \eqref{pit} that $\hat{E}^v$ is in fact a Lubin-Tate group
over ${\cal O}_v$ for the local parameter $\pi.$

\begin{lem}\label{qot}
Let $A_{\lambda}(t)$ be the $t$-expansion of $\Psi_{\lambda}(P).$
Then $A_{\lambda}(t)$ belongs to $1+q{\cal O}_v[[t]].$
\end{lem}

\begin{cor}\label{clt}
Defining
\begin{equation}\label{cln}
C_{\lambda}(t)\,=\,(Nv)^{-1}\log\,A_{\lambda}(t),
\end{equation}
we have $q\,C_{\lambda}(t)$ belongs to ${\cal O}_v[[t]].$
\end{cor}
\proof
Let $B_{\lambda}(t)$ be the $t$-expansion of ${\frak
  R}_{\lambda}(P)$.
A standard argument (see the proof of Lemma 8 in \cite{LMS}),
based on the explicit expression \eqref{r1}, shows that
$B_{\lambda}(t)$
is a unit in ${\cal O}_v[[t]].$ It follows that
$B_{\lambda}([\pi](t))$
is also a unit in ${\cal O}_v[[t]],$ whence we conclude from
\eqref{psil}
that $A_{\lambda}(t)$ belongs to ${\cal O}_v[[t]].$ Moreover,
writing
$B_{\lambda}(t)=\,\overset{\infty}{\underset{n=0}{\Sigma}}\,b_nt^n,$
we deduce from \eqref{pit} that
\begin{equation}\label{blp}
B_{\lambda}([\pi](t))\,=\,\overset{\infty}{\underset{n=0}
{\Sigma}}\,b_n([\pi](t))^n \equiv \overset{\infty}{\underset{n=0}
{\Sigma}}\,b_nt^{nNv}\mod q.
\end{equation}
On the other hand, as $b_n^{Nv}\equiv b_n \mod q$, one has
\begin{equation}\label{bln}
B_{\lambda}(t)^{Nv} \equiv \overset{\infty}{\underset{n=0}{\Sigma}}
\,b_nt^{nNv}\mod q.
\end{equation}
It follows immediately from \eqref{blp} and \eqref{bln} that
$A_{\lambda}(t) \equiv 1 \mod q,$ completing the proof of
Lemma \ref{qot}. The corollary is immediate since $\log\,A_{\lambda}(t)$
then belongs to $q{\cal O}_v[[t]]$, since $q$ is odd.
\qed

Define
\begin{equation}\label{26}
C_{\lambda}^*(t)=\,q C_{\lambda}(t).
\end{equation}
By Corollary \eqref{clt}, $C_{\lambda}^*(t)$ belongs to
${\cal O}_v[[t]],$ and, by \eqref{psilp},
\begin{equation}\label{27}
\underset{R \in E_q}{\sum}\, C_{\lambda}^*\left(t[+]t(R)\right)=0,
\end{equation}
where $[+]$ denotes the formal group law on $\hat{E}^q.$ Here
$t(R)$ runs over the $q$-division points on $\hat{E}^q$ as $R$
runs over $E_q.$ Katz's argument  applies to any power series in ${\cal O}_v[[t]]$ satisfying
\eqref{27}.

\begin{prop}\label{2.5}(Katz [K])
Let $g(t)$ be any power series in ${\cal O}_v[[t]]$ satisfying
\begin{equation}\label{28}
\underset{R \in E_q}{\Sigma}\, g\left(t[+]t(R)\right) =0.
\end{equation}
Then, for all integers $n \geq 1, $ we have
$$
\left(\ddz\right)^n g(t)\, \in q^{[\frac{nq}{q^2-1}]}\,{\cal O}_v[[t]].
$$
\end{prop}
We now briefly describe Katz's proof. It is convenient to replace $\hat{E}^q$ by an isomorphic Lubin-Tate group.
Let ${\cal E}$ be the Lubin-Tate group over ${\cal O}_v$ attached to the local parameter $\pi$, and satisfying
\begin{equation}\label{29}
[\pi](w)=\pi w + w^{Nv},
\end{equation}
where $[\pi](w)$ now denotes the endomorphism of ${\cal E}$
defined by $\pi.$ By Lubin-Tate theory, there exists an ${\cal O}_v$-isomorphism
\begin{equation}\label{30}
\eta\,: \,{\cal E} \simeq \hat{E}^q,
\end{equation}
which is given by a formal power series $t=\eta(w)$ in
in ${\cal O}_v[[w]].$ Defining $h(w)=g(\eta(w)),$ we then have
\begin{equation}\label{31}
\underset{u \in {\cal E}_q}{\Sigma}\, h(w[+] u)=0,
\end{equation}
where $[+]$ also denotes the formal group law on ${\cal E}$, and ${\cal E}_q$ is the group of $q$-division points on ${\cal E}.$ The isomorphism \eqref{pz} enables us to write
$z=\varepsilon(t),$ where $\varepsilon (t)$ is a power series in $K_v[[t]].$ It is then easy to see that
\begin{equation}\label{32}
z= \nu(w), ~{\rm where}~\nu(w)=\varepsilon(\eta(w))
\end{equation}
is the logarithm map of ${\cal E}$, and that it suffices to show that, for all integers $n \geq 1$,
\begin{equation}\label{33}
\left(\ddz\right)^n\,h(w)\,\in q^{[\frac{nq}{q^2-1}]}\,{\cal O}_v[[w]].
\end{equation}
Note that the operator
\begin{equation}\label{34}
\ddz\,=\,\frac{1}{\nu '(w)}\frac{d}{dw}
\end{equation}
maps ${\cal O}_v[[w]]$ into itself, since $\nu '(w)$ is a unit in ${\cal O}_v[[w]]$ by a well known property of formal groups.

\medskip

For each $r(w)$ in ${\cal O}_v[[w]]$, we define $D_n\,r(w)~(n \geq 0)$ in
${\cal O}_v[[w]]$ by the expansion
\begin{equation}\label{35}
r(w[+]u)\,=\,\overset{\infty}{\underset{n=0}{\Sigma}}\, D_n\,r(w)u^n,
\end{equation}
where $u$ is an independent variable.

\begin{lem}\label{2.6}
For $n=0,\,1,\cdots, Nv-1,$ we have
\begin{equation}\label{36}
D_n\,r(w)\,=\, \frac{1}{n!} \left(\ddz\right)^n r(w).
\end{equation}
\end{lem}
\proof
Since $r(w[+]u_1[+]u_2) = r(w[+]u_2[+]u_1)$,  we obtain
\begin{equation}\label{37}
D_{n_1}(D_{n_2}\,r(w))\,=\, D_{n_2}\left(D_{n_1}\,r(w)\right)~~~~
(n_1,\,n_2\geq 0).
\end{equation}
Also, as $\nu(w[+]u)=\nu(w)+ \nu(u),$ we have
$$
\frac{\partial}{\partial u} \, (w[+]u)\,=\, \nu'(u)/\nu'(w[+]u),
$$
and hence
$$
\frac{\partial}{\partial u} \, r(w[+]u)\,=\,
r'(w[+]u)\nu'(u)/\nu'(w[+]u).
$$
Putting $u=0$ in this equation, and noting that $\nu'(0)=1$, it follows that $D_1r(w) = r'(w)/\nu'(w).$
Thus the above equation can be rewritten as
\begin{equation}\label{38}
\frac{\partial}{\partial u}\,r(w[+]u)\,=\,\nu'(u)(D_1\,r)(w[+]u).
\end{equation}
Recalling \eqref{37}, this then gives the identity
\begin{equation}\label{39}
\overset{\infty}{\underset{n=1}{\Sigma}}\, n D_n \,r(w) u^{n-1}=\nu '(u)
\overset{\infty}{\underset{n=0}{\Sigma}}\, D_1(D_n\,r((w))u^n.
\end{equation}
Also, as $\nu([\pi](w))=\pi\nu(w),$ one easily deduces from \eqref{29} that
\begin{equation}\label{40}
\nu '(w) \equiv 1\!\mod w^{Nv-1}.
\end{equation}
Combining \eqref{39} and \eqref{40}, we immediately obtain $D_n\,r=\,D_1(D_{n-1} r)/n$ for $n=1,\cdots, Nv-1,$
and the assertion of the lemma follows by induction on $n$.
\qed

\medskip

Since $\ord_q((Nv-1)!)=q-1,$ it follows immediately from Lemma
\ref{2.6} that, for each $r(w)$ in ${\cal O}_v[[w]],~ \left(\ddz\right)^{Nv-1}\,r(w)$ belongs to $q^{q-1}{\cal O}_v[[w]].$ The next lemma establishes a stronger result, provided $r(w)$ satisfies
\begin{equation}\label{41}
\underset{u \in {\cal E}_q}{\Sigma}\, r(w [+] u)=0.
\end{equation}
We note that if $r(w)$ satisfies \eqref{41}, so also does $\left( \ddz\right)^nr(w)$ for all non-negative integers $n$, since
$\nu'(w[+]u)=\nu'(w)$.

\begin{lem}\label{2.7}
Let $r(w)$ be any power series in ${\cal O}_v[[w]]$ satisfying \eqref{41}.
Then
\begin{equation}\label{42}
\left( \ddz\right)^{Nv-1}\,r(w)\in q^q\,{\cal O}_v[[w]].
\end{equation}
\end{lem}
\proof
Combining \eqref{35} and \eqref{41} gives
\begin{equation}\label{43}
\overset{\infty}{\underset{n=0}{\Sigma}}
\left(\underset{\eta\in{\cal E}_q}{\Sigma}\,\eta^n\right)\, D_nr(w)=0.
\end{equation}
By \eqref{29}, we see that the non-zero elements of ${\cal E}_q$ are given by the $\alpha \zeta,$ where $\zeta$ runs over the $(Nv-1)$-th roots of unity, and $\alpha^{Nv-1}=-\pi.$ Thus the $n$-th term in
\eqref{43} is zero unless $Nv-1$ divides n, and, when $Nv-1$ does divide $n$, we have
$$
\underset{\eta\in{\cal E}_q}{\Sigma}\,\eta^n \,=\, (Nv-1)\alpha^n
~~~~~(n > 0).
$$
Hence \eqref{43} can be written as
\begin{equation}\label{44}
q^2r(w) + \overset{\infty}{\underset{m=1}{\Sigma}}\, (Nv-1)(-\pi)^m\,
D_{m(Nv-1)}\,r(w)=0.
\end{equation}
But this last equation clearly implies that
$$
D_{Nv-1}\,r(w)\in q{\cal O}_v[[w]],
$$
and the assertion of the lemma now follows from Lemma \ref{2.6} and the remarks made immediately before \eqref{41}.
\qed

\begin{lem}\label{2.8}
Let $r(w)$ be any element of ${\cal O}_v[[w]]$ satisfying \eqref{41}.
Then, for all integers $n \geq 0,$ we have
\begin{equation}\label{45}
\left(\ddz\right)^n\,r(w) \in q^{[\frac{nq}{Nv-1}]}\,{\cal O}_v[[w]].
\end{equation}
\end{lem}
\proof
Assume $n \geq 1,$ and write
$$
n =(Nv-1)b +a,
$$
where $b\geq 0$ and $0 \leq a < Nv-1.$ Now
$$
\left[\frac{nq}{Nv-1}\right]\,=\,\left[ \frac{aq}{Nv-1}\right] +bq,
$$
and
$$
\left[\frac{aq}{Nv-1}\right]\,=\, \ord_q(a!).
$$
The assertion of the lemma is now plain from \eqref{36} and
\eqref{42}.
\qed

\medskip

We can now complete the proof of Theorem \ref{kat}. Proposition \ref{2.5}
follows immediately from applying Lemma \ref{2.8} to the function
$h(w)$. In turn, applying Proposition \ref{2.5} to the power series
$C_{\lambda}^*(t)$ given by \eqref{26}, and recalling that $(f,v)=1$, we conclude from Proposition \ref{ddz} that, for all integers $n \geq 1$,
$$
\ord_v\left((N\lambda-\psi_E((\lambda))^n)\left(1-\frac{\psi_E(v)^n}{Nv}\right)L_n\right) \geq \left[\frac{nq}{Nv-1}\right]-1.
$$
Assuming that $n \geq 3,$ it is clear that $1-\left(\psi_E(v)^n/Nv\right)$ is a $v$-adic unit. Now choose $\lambda$ to be any element of ${\cal O}_K$ such that $(\lambda,6)=1,$ $\lambda \equiv 1 \!\mod {\frak f},$ and $\lambda \mod v$ is a generator of $\left({\cal O}_K/v {\cal O}_K\right)^{\times}.$ As $\lambda \equiv 1 \!\mod {\frak f},$ $\psi_E((\lambda))=\lambda,$ and thus
\begin{equation}\label{46}
N\lambda-\psi_E((\lambda))^n=\lambda(\bar \lambda-\lambda^{n-1}) \equiv
\lambda(\lambda^q-\lambda^{n-1}) \! \mod \,v.
\end{equation}
Since $\lambda$ mod $v$ is a generator of $\left({\cal O}_K/v{\cal O}_K\right)^{\times},$ it
follows that \eqref{46} is prime to $q$ provided $n$ is not congruent to $1+q$ mod $Nv-1.$ This completes the proof of Theorem \ref{kat}.
\qed

\section{Application to the Tate-Shafarevich group}

In this section, we combine Theorem \ref{kat} with Iwasawa theory to prove
Theorem \ref{1.1}. Throughout  the section, $c_1,\,c_2,\cdots$ will
denote positive integers which depend only on the coefficients
of the equation \eqref{eqe}, and which
could be made explicit if desired. Also, $J$ will denote the set
of all prime numbers $q$ such that $q$ is inert in $K$, and $E$
has good reduction at $q{\cal O}_K.$

\medskip

For the moment, let $p$ be any odd prime number. Put
\begin{equation}\label{47}
{\cal P}(p) =\underset{\underset{q\leq p}{q\in J}}{\prod}\,q^{\nu_q},
~~{\rm where}~
\nu_q=\left[\frac{pq}{q^2-1}\right]-1.
\end{equation}

\begin{lem}
For each odd prime $p$, we have
\begin{equation}\label{47'}
{\cal P}(p) \geq p^{p/2}/c_1^p.
\end{equation}
\end{lem}
\proof
Clearly
$$
\nu_q \geq \left[\frac{p}{q}\right]-1 \geq \frac{p}{q}-2,
$$
and so
\begin{equation}\label{48}
\log\,{\cal P}(p) \geq p \underset{\underset{q \leq p}{q \in J}}{\sum}\, \frac{\log\,q}{q}- 2\underset{\underset{q \leq p}{q \in J}}{\sum}\, \log\, q.
\end{equation}
By a weak form of the prime number theorem, we have
\begin{equation}\label{49}
\underset{\underset{q \leq p}{q \in J}}{\sum}\,\log \,q
< \underset{\underset{q {\rm prime}}{q \leq p}}{\sum}\,
\log\, q \leq c_2 p.
\end{equation}
Let $\chi$ be the Dirichlet character corresponding to the extension $K/\Q$, and
let $t$ be its conductor. Define
\begin{equation}
{\cal M}(p)\,=\, \underset{\underset{q \leq p}{\chi(q)=-1}}{\sum}\,
\frac{\log\,q}{q},
\end{equation}
where the sum is taken over all primes $q \leq p$ with $\chi(q)=-1.$ Plainly
\begin{equation}\label{50}
0 \leq {\cal M}(p) - \underset{\underset{q \leq p}{q \in J}}{\sum}\,
\frac{\log q}{q} \leq c_3.
\end{equation}
Combining \eqref{48}, \eqref{49}, and \eqref{50}, we obtain
\begin{equation}\label{51}
\log\,{\cal P}(p) \geq p{\cal M}(p) - c_4p.
\end{equation}
Now a well known equivalent form of Dirichlet's  theorem on primes in arithmetic progressions asserts that,
for any real $x \geq 2$, and each integer $q$ with $(a, t)=1,$
we have
\begin{equation}\label{52}
\left|\underset{\underset{q \leq x}{q \equiv a \!\!\mod t}}{\sum}\,
\frac{\log\,q}{q}\, - \frac{1}{e(t)}\, \log \, x \right| \leq c_5,
\end{equation}
where the sum on the left is over all prime numbers $q \leq  x$ with $q \equiv
a \mod t,$  and $e(t)$ denotes the order of the group of units of  $\Z/t\Z$.  Now precisely
half of the classes $a$ mod $t$ with $(a,t)=1$
satisfy $\chi(a)=-1.$ Taking $x=p$, we conclude  from \eqref{52} that
\begin{equation}\label{53}
\mid {\cal M}(p) - \frac{\log\,p}{2}\mid \leq c_6.
\end{equation}
The assertion of the lemma now follows immediately on
combining \eqref{51} and \eqref{53}.
\qed

\medskip

The proof of Theorem \ref{1.1} now proceeds exactly as in the proof of Theorem \ref{2.8}
of \cite{CZS}, except we exploit Katz's divisibility assertion \eqref{super}. Thus we now take $p$ to be any odd prime with $p{\cal O}_K={\frak p}{\frak p}^*$ and ${\frak p}\neq {\frak p}^*,$  such that $E$ has good reduction at both ${\frak p}$ and
${\frak p}^*.$ Then there exists a positive rational integer $c_7$, depending only on the equation \eqref{eqe}, such that
$$
c_7^p(p-1)!\,\Omega_{\infty}^{-p}L_{\frak f}(\bar{\psi}_E^p,p)
$$
is an algebraic integer in $K$ (see the proof of Lemma \ref{2.8} of \cite{CZS}). This algebraic integer is non-zero because the Euler product for $ L_{\frak f}(\bar{\psi}_E^p,s)$  converges at $s=p.$ Moreover, by Theorem \ref{kat}, it is divisible by the rational integer ${\cal P}(p).$ Thus
\begin{equation}\label{54}
\xi(p)\,=\, c_7^p(p-1)!\,\Omega_{\infty}^{-p}\,L_{\frak f}(\bar{\psi}_E^p,p)/{\cal P}(p)
\end{equation}
remains a non-zero algebraic integer in $K$. In view of Lemma 3.1 and \cite[Lemma 2.9]{CZS}, we have
\begin{equation}\label{55}
\mid\xi(p)\mid \, \leq \, c_8^p\, p^{p/2}.
\end{equation}
Thus, using the product formula for $K$, we conclude that
\begin{equation}\label{56}
 \mid\xi(p)\mid_{\frak p} \, . \,  \mid\xi(p)\mid_{\frak p^*} \, \geq \, c_8^{-2p}p^{-p} \, \geq \, p^{-(1+\epsilon)p};
\end{equation}
here the last inequality holds for any $\epsilon > 0$, provided $p$ is sufficiently large. On the other hand, applying the main conjectures  of Iwasawa theory (which are in fact proven theorems) for $E$ over the unique $\Z_p$-extensions of $K$ which are unramified outside $\frak p$ and $\frak p^*$, respectively, it follows that (see the proof of Theorem 2.1 in \cite{CZS})
\begin{equation}\label{57}
 \mid\xi(p)\mid_{\frak p} \, . \,  \mid\xi(p)\mid_{\frak p^*} \, \leq \, p^{-2(n_{E/K}+t_{\frak p})},
\end{equation}
where $n_{E/K}$ is the ${\cal O}_K$-rank of $E(K)$, and $t_{\frak p}$
is the $\Z_p$-corank of $\sha(E/K)(\frak p)$
(which can be shown to be equal to the $\Z_p$-corank of
$\sha(E/K)(\frak p^*)$).
Thus, combining \eqref{56} and \eqref{57}, we have proven the following stronger form of Theorem 2.8
of \cite{CZS} :-

\begin{thm}\label{58}
Let $\epsilon$  be any positive number. Then, for all sufficiently
large odd primes
 $p$ such that  $p{\cal O}_K={\frak p}{\frak p}^*$, $t_{\frak p}$
is bounded above by $(1/2 + \epsilon)p - n_{E/K}.$
\end{thm}
Finally, we note  that Theorem \ref{1.1} is an immediate
corollary of this result, because, when $E$ is defined over $\Q$
we have $n_{E/K} = g_{E/\Q}$ and $t_{\frak p} = t_{E/\Q,p}$. This completes the proof of Theorem \ref{1.1}.
\qed

\section{Computations for $y^2=x^3-Dx$.}

We explain in this section the improvements in the computational technique
of \cite{CZS}, which enables us to prove Theorem \ref{1.2} and
Theorem\ref{ei}.
As in \cite{CZS}, we consider the family of curves
$$
E\,:\, y^2=x^3-Dx,
$$
where $D$ is a fourth power free non-zero rational
integer. For these curves, $K=\Q(i)$ and the isomorphism from
$\Z[i]$ to $\End_K(E)$ is given by mapping $i$ to the endomorphism which sends $(x,y)$ to $(-x,iy).$

\medskip

As earlier, let ${\frak f}$ be the conductor of the Gr\"ossencharacter
$\psi_E$, and fix some $f$ in ${\cal O}_K$ such that ${\frak f}=f{\cal
  O}_K$. The explicit value of ${\frak f}$ is well-known, and is given
by Lemma 3.2 of \cite{CZS}. In particular, ${\frak f}$ is always
divisible by $(1+i){\cal O}_K$, and we define
\begin{equation}\label{60}
f_1=f/(1+i),~~{\frak f}_1=f_1{\cal O}_K.
\end{equation}
Recall that $G$ denotes the Galois group of $K(E_{\frak f})$ over $K$.
If $\frak g$ is an integral ideal of $K$, we denote the order of the multiplicative group
of units of ${\cal O}_K/\frak g$ by $\phi(\frak g)$. Also the symbol $\ominus$ will denote
subtraction in the group law of $E$.

\begin{lem}
Assume that $D$ is divisible by an odd prime. Let $H_{{\frak f}_1}$  be the ray class field of $K$  modulo ${\frak f}_1$.  Then  $H_{{\frak f}_1} = K(u)$ , where $u$ is the $x$-coordinate of any primitive ${\frak f}_1$-division point
on $E$. Moreover, $[K(E_{\frak f}):H_{{\frak f}_1}]=2$, and
\begin{equation}\label{61}
\sigma(Q)=\ominus Q~~{\rm {for~all}}~ Q~{\rm in}~ E_{{\frak f}_1},
\end{equation}
where $\sigma$ denotes the non-trivial element  of  the Galois group of $K(E_{\frak f})$ over $H_{{\frak f}_1}$.
\end{lem}
\proof
Since $D$ is odd, ${\frak f}_1$ is divisible by a prime of $K$
distinct from $(1+i){\cal O}_K$, and thus
if $\zeta$ is a root of unity in $K$ with $\zeta \equiv 1 \mod\!{\frak
  f}_1$, then we must have $\zeta=1.$ It follows that
$[H_{{\frak f}_1}:K]=\phi({{\frak f}_1})/4.$
On the other hand, since $\frak f$ is the conductor of $\psi_E$,
$K(E_{\frak f})$ coincides with the ray class field of $K$
modulo ${\frak f}$ (see \cite[Lemma 7]{LMS}), and so
$[K(E_{\frak f}):K]=\phi(\frak f)/4$,  whence
\begin{equation}\label{62}
[K(E_{\frak f}):H_{{\frak f}_1}]=2.
\end{equation}
Let$(u,v)$ be any
primitive ${{\frak f}_1}$-division point on $E$. By the classical
theory of complex multiplication, $H_{{\frak f}_1}=K(u^2).$
Let $\sigma$ denote the non-trivial element of $\Gal(K(E_{\frak
  f})/H_{{\frak f}_1}).$
Since $H_{{\frak f}_1}=K(u^2),$ we must have $\sigma\,u=\pm u.$ But
$\sigma\,u=-u$ is impossible, since it would imply that $\sigma\,v=\pm
iv,$ which would in turn imply that $\sigma^2$ is not the identity
element of $G$. Hence
$\sigma\,u=u.$ The following argument shows that we cannot have
$\sigma v=v$. Assume indeed that $\sigma v=v$, so that
\begin{equation}\label{63}
H_{{\frak f}_1}=K(E_{{\frak f}_1}).
\end{equation}
Take $\alpha$ to be any element of ${\cal O}_K$ so that $\alpha \equiv
1 \mod \! {\frak f}_1,$ and put ${\frak a}=\alpha{\cal O}_K.$ Note
that ${\frak a}$ is prime to ${\frak f}$ because ${\frak f}$ and
${\frak f}_1$ have the same prime factors. Since the abelian extension
$H_{{\frak f}_1}/K$ has conductor ${\frak f}_1$, the Artin symbol
  $\tau_{\frak a}$ of ${\frak a}$ for the extension $H_{{\frak f}_1}/K$
is equal to 1. But by \eqref{gr1}, we have
\begin{equation}\label{64}
\tau_{\frak a}(Q)=\psi_E({\frak a})(Q)~~{\rm for~all}~  Q ~{\rm
  in}~E_{{\frak f}_1}.
\end{equation}
Hence we must have $\psi_E({\frak a})\equiv 1 \!\mod {\frak f}_1.$ But
${\frak a}=\psi_E({\frak a}){\cal O}_K,$ and so $\psi_E({\frak
  a})=\zeta\alpha$ where $\zeta$ is a root of unity in $K$.  As
$\alpha \equiv 1\!\mod {\frak f}_1,$ it follows that
$\psi_E((\alpha))=\alpha$ for all $\alpha$ in ${\cal O}_K$ with
$\alpha \equiv 1 \! \mod {\frak f}_1.$
This in turn implies
that the conductor ${\frak f}$ of $\psi_E$ must divide ${\frak f}_1$,
which is a contradiction. Hence we must have
$$
\sigma(u,v)=(u,-v)=\ominus(u,v),
$$
which proves \eqref{61}.  This completes the proof.
\qed

\medskip

The Weierstrass differential equation associated to $E$ is
\begin{equation}\label{65}
{\wp}'(z,{\cal L})^2=4\wp(z,{\cal L})^3-4 D\wp(z,{\cal L}).
\end{equation}
In general, we write ${\wp}^{(n)}(z,{\cal L})$ for the $n$-th
derivative of $\wp(z,{\cal L})$ with respect to $z$. For all integers
$n \geq 0$, one has
\begin{equation}\label{66}
\wp^{(2n+1)}(z,{\cal L})=B_n(\wp(z,{\cal L}))\wp'(z,{\cal L}),
\end{equation}
where $B_n(X)$ is a polynomial of degree $n$ in $\Z[X].$

\begin{cor}\label{4.2}
Assume that $D$ is divisible by an odd prime. Then, for all integers $n
\geq 0$,
\begin{equation}\label{67}
{\rm Tr}_{K(E_{\frak
    f})/K}\,\left(\wp^{(2n+1)}
(\Omega_{\infty}/f_1,{\cal
      L})\right)=0.
\end{equation}
\end{cor}
\proof
Writing ${\rm Tr}_{K(E_{\frak f})/H_{{\frak f}_1}}$ for the trace map
from $K(E_{\frak f})$ to $H_{{\frak f}_1}$, it is clear from
\eqref{61} and \eqref{66} that
\begin{equation}\label{68}
{\rm Tr}_{K(E_{\frak f})/H_{{\frak f}_1}}\left
(\wp^{(2n+1)}(\Omega_{\infty}/f_1,{\cal
      L})\right)=0.
\end{equation}
The assertion \eqref{67} follows immediately, completing the proof.
\qed

\medskip

We next introduce the formal expressions
\begin{equation}\label{69}
W(z)=\wp(z,{\cal L})^{1/2},~~~V(z)=\left(\wp(z,{\cal L})^2-D
\right)^{1/2}.
\end{equation}
Noting that the differential equation \eqref{65} can be written as
\begin{equation}\label{dif}
\wp'(z,{\cal L})=2 W(z) V(z),
\end{equation}
one immediately obtains
\begin{equation}\label{70}
W'(z)=V(z),~~V'(z)=2 W(z)^3.
\end{equation}
The following lemma is then clear by induction on $n$.

\begin{lem}\label{4.3}
For all integers $n \geq 0$, we have
$$
V^{(2n+1)}(z)\,=\, A_n\left(W(z)\right),
$$
where $A_n(X)$ is a polynomial of degree $2n+3$ in $\Z[X].$
\end{lem}

\medskip

\noindent  Since
$$
\wp\left((1+i)z,{\cal L}\right)=\frac{\wp(z,{\cal
    L})^2-D}{2i\wp(z,{\cal L})},
$$
it follows that
\begin{equation}\label{71}
\wp(z,{\cal L})=i\,\wp\left((1+i)z,{\cal L}\right)
\pm  i\,V\left((1+i)z\right),
\end{equation}
Define
\begin{equation}\label{ro}
\rho=W\left(\frac{\Omega_{\infty}}{f_1}\right),
\end{equation}
so that $\rho^2=\wp(\Omega_{\infty}/f_1 , \cal L ).$

\begin{lem}\label{ro1}
Assume that D is divisible by an odd prime.
Then $K(\rho) = H_{{\frak f}_1}$, where  $H_{{\frak f}_1}$ is the ray class field of $K$ mod $\frak f_1$.
In particular,  $[K(\rho):K] = \phi(\frak f)/8$.
\end{lem}
\proof
Putting $ z=\Omega_{\infty}/f$ in \eqref{71}, we conclude that $\gamma = V(\Omega_{\infty}/f_1)$ belongs to
$F = K(E_{\frak f})$.   Taking
$ z=\Omega_{\infty}/f_1$ in \eqref{dif},  it then follows that
$\rho$  must also belong to $F = K(E_{\frak f})$. Moreover, as
$F=K( \wp(\Omega_{\infty}/f , \cal L))$, we see from \eqref{71}
that $\gamma$ does not belong to  $H_{{\frak f}_1}$, whereas
its square plainly does.  Hence,  writing $\sigma$  for the
non-trivial element of  $Gal(F/H_{{\frak f}_1})$, it follows
that $ \sigma(\gamma) = -\gamma$.  We know from \eqref{61} that we also have
$\sigma(v) = - v$, when $ v= \wp'(\Omega_{\infty}/f_1 , \cal L)$.
It then follows from \eqref{dif} with  $ z=\Omega_{\infty}/f_1$
that $\sigma$ must fix $\rho$, completing the proof.
\qed

\medskip

We continue to assume that $D$ is divisible by an odd prime.  Then, after differentiating \eqref{71}, and using \eqref{67} and the lemma just proven,  we conclude that, for all integers $n \geq 0$,

\begin{equation}\label{72}
{\rm Tr}_{K(E_{\frak
    f})/K}(\wp^{(2n+1)}(\Omega_{\infty}/f,{\cal L}))
=\,
\pm (1+i)^{2n+3}\,{\rm Tr}_{H_{\frak
    f_1}/K}\left(A_n(\rho)\right).
\end{equation}
We use the right hand side of this formula to compute the left hand side in our
numerical calculations. The great  advantage from a numerical point of view  is that
the equation  \eqref{72} allows us to avoid knowing the minimal polynomials of both
$\wp(\Omega_{\infty}/f,{\cal L})$ and $\wp'(\Omega_{\infty}/f,{\cal
  L})$  over $K$, and only requires knowledge of the minimal polynomial of
$\rho$ as well as the coefficients of $A_n(X).$ Indeed, as we shall explain in the next paragraph, we have been able to calculate
the minimal polynomial of $\rho$ over $K$ for each of  the 6 curves
discussed in the Introduction.

\medskip

Let $H(X)$ be the minimal polynomial of $\rho$ over $K$. By Lemma
\ref{ro1}, $H(X)$ has degree $\phi({\frak f})/8.$ We succeeded in
computing $H(X)$ when
\begin{equation}\label{76}
D\,=\,-14,\,17,\,-33,\,-34,\,-39, \,82,
\end{equation}
and the respective degrees of $H(X)$ are 192, 128, 960,
1024, 576, 64,000. By far the most difficult case was $D=82,$
and it was  only found with the help of the
MAGMA programmes. It would take too much space to give these polynomials
explicitly here, although we would be happy to provide an electronic
file containing them for any interested reader. We now list the
coefficient of largest absolute value for each of the polynomials $H(X)$,
together with the power of $X$ where it occurs:-
\medskip

\noindent $D=-14$; coefficient of $X^{40}$:
\begin{align*}-17505698603459355436042213669487147582723
059629657863703494656.\end{align*}

\noindent $D=17$; coefficient of $X^{18}$:
\begin{align*}-323854307090694728597766056638496367408758560.\end{align*}

\medskip

\noindent $D=-33$; coefficient of $X^{107}$:
\begin{align*}&
1113692168077398337776583577672684033100787524823564456336793469544236\\
&9882566874514629785310533157148094245351329286533498412231880301090833\\
&6242433502572436644559991342707201233037605661492070818654076942778326\\
&7556086484634044977612377544550389135642352731231407460900676813436679\\
&7044082034896262798652340498136472798865345948403494563087565937673700\\
&28550083854901398254552448895886221992192.
\end{align*}

\medskip

\noindent $D=-34$; coefficient of $X^{104}$:
\begin{align*}&
7720430792855503034710415580894323862364569575047279968634760866523454\\&
7847422997182187637013858030921907403287538876053146805325940886391064\\&
7892199987049394255653536333933477502477411283490330578317508812402336\\&
2555487142551124594924227469513073746361763633056985881176027866067892\\&
1523140860877514259352431551963116165178346512108029927115844217508794\\&
334700569430116613294495880552612844403970375491308137706684416.
\end{align*}

\medskip

\noindent $D=-39$; coefficient of $X^{52}$:
\begin{align*}&
-355740716783970448334076807956400505497970442768469601110698516004184\\&
8127114241683240295186562457896490055695622807539265131729992420618358\\&
6571720641221008764113343860595870542652478897038226690477175270591499\\&
924870328538293899715218933146816(1+i).
\end{align*}

\medskip

\noindent $D=82$; coefficient of $X^{376}$:
\begin{align*}&
1626770346973377213003989713931931718932500126487797045408201501253939\\&
6731694266094747666615530988064086303049611694143030937780501940450482\\&
6148923158148399191147813882230238624119655402654012877646848571872713\\&
2822553302999268257535786699238655250299150648286294022464086021755163\\&
9738764276240402326739779256978620234034249213967253713905896267785931\\&
1389272946311927380327203467553385171814069136924311033952190656371048\\&
8401256493656625338991614392352334979047794680542663138130567100068261\\&
4730080755894768964763705148318306057256912022688939397372813267988367\\&
6652353478038859315100790273860531016767864764445838519619890332349983\\&
1711986989627104655179084273340230013821136292556801493781093121612681\\&
4363981856242791856465666624336990116907158556227793717809986779598196\\&
1315873561402664253429232577346034909009263102014156414450768913713516\\&
8222227131939890289476123528846055893652489735995479793826504226262779\\&
2807348394599705036406983089127665543492586507207678193393235364804223\\&
6462816194025382172710423211077822337257079817790280658424828864503574\\&
7985383293171188977385859941751648239156399693931547246990152507384535\\&
0320264105620916081385901637880964310002635955229618201426636587238896\\&
1043298009429916985454588056855950483173445938068590086144429016871410\\&
\end{align*}
\newpage
\begin{align*}&
2301521229865047393894107278179025053626765685604299939779709329204331\\&
6117427016830769698566119962791577731251538197220739446420738649715257\\&
9550649615272922868402516863443993334221538823660227843895865513813349\\&
0677958099100627837595866989226624619466138791950725644948832889480533\\&
6196454453562854270955682967252379054309054191951467056681904480077015\\&
4504643277448756066354911910313550003764314073875751467126140895984307\\&
9734334324608414916401831949815825979625541958282127234282122960726210\\&
1911031290501410725395455324594186221966782329174986332874204939357635\\&
8485688197981820411047322520602151422097104881189661486945477984943401\\&
2544237753866264395719050412246974517882381632563194789563672511042082\\&
9052803406635594768668706675481156888453870362673285085628168472992950\\&
1901225272744412285530261592450972025459750324518063311741511320093528\\&
2619130188893743499867860544264225634098941938446218692450707841031742\\&
5542145118294934165120179815306657244830115173958401214860691464241382\\&
5695303908587970350882714782329517816083108740133499352732737947890551\\&
1840108553923823745181741491616613465575991279099569495745640109613645\\&
5013082046639240884294500481314527388524452661343943844013075136023422\\&
6584580874160255211038190608666792275054918912669599928553999654855202\\&
7886831115375661963183338778623436402481583278842727187320847197910254\\&
6899871231485102958986194253260192603773895356652620176464371298043769\\&
1367330290209445179016520776412957412615689757650463346453062234881546\\&
7129189078717972033343356971946046854148898545124234537512138427279216\\&
5082778329041394790591607367236385469309097587612551453985632669821803\\&
0228499100164222610066790856089524618448511113564527894301127242918770\\&
2759120062569482117569076417828776520673334475477832948015993478979289\\&
3767008605869336671146314674407916016696140772323044613933507846739257\\&
3030870182550866734231510577695191594238021017596602280094728192.
\end{align*}

\medskip

We now discuss briefly the computation of the relevant $L$-values.
Let $\Omega_{\infty}^+$ be the smallest positive real period of $E$, so
that
$$
\Omega_{\infty}=\Omega/D^{1/4}~{\rm if}~ D>0,~~\Omega_{\infty}^+=
\Omega/(-D/4)^{1/4}~{\rm if}~ D<0,
$$
where $\Omega=2.622058...$ is the least real non-zero period of
$y^2=x^3-x.$ Thus $\Omega_{\infty}^+=\alpha(E)\Omega_{\infty},$
where $\alpha(E)$ is 1 or $(1+i)$ according as $D>0$ or $D<0.$
As in \cite{CZS}, we write
\begin{equation}\label{77}
c_n^+(E)=(\Omega_{\infty}^+)^{-n}\, L(\bar{\psi_E}^n,n)~~~(n \geq 1),
\end{equation}
so that $c_n^+(E)$ is in $\Q.$

\medskip

Now take $p$ to be any prime with $p \equiv 1$ mod 4 and $(p,D)=1.$
Using the full force of the main conjectures of Iwasawa theory, it is
proven in \cite[Theorem 2.2]{CZS} that (i) $\ord_p(c_p^+(E))\geq g_{E/\Q},$
(ii) if $\ord_p(c_p^+(E))=g_{E/\Q},$ then $\sha(E/\Q)(p)=0,$ and (iii)
if $\ord_p(c_p^+(E))\leq g_{E/\Q}+1,$ then $\sha(E/\Q)(p)$ is finite. We
calculate $\ord_p(c_p^+(E))$ as follows. Putting $n=p$ in \eqref{ln}, one has
$$
c_p^+(E)\,=\,\frac{1}{\alpha(E)^p(p-1)!}\,L_p,
$$
and, by a classical formula,
$$
L_p= -f^{-p}{\rm Tr}_{K(E_{\frak f})/K}\,
\left(\wp^{(p-2)}\left(\frac{\Omega_{\infty}}{f},{\cal L}\right)\right).
$$
Combining this with \eqref{72} for $n=(p-3)/2$, we obtain
\begin{equation}\label{78}
c_p^+(E)=\pm \frac{1}{\beta(E)^p(p-1)!}\,{\rm Tr}_{H_{{\frak f}_1}/K}\,
\left(A_{\frac{p-3}{2}}(\rho)\right),
\end{equation}
where
\begin{equation}\label{79}
\beta(E)=f{\alpha(E)}/(1+i).
\end{equation}
Of course, $c_p^+(E)$ is obviously positive, and we choose the
sign in \eqref{78} accordingly. This formula \eqref{78} is the key
one for the computations. Indeed, the expression on the right of
\eqref{78} can readily be calculated when we know explicitly the minimal
polynomial of $\rho$ over $K$.

\medskip

For the six values of $D$ given by \eqref{76}, we make use of the
following respective values of $f$:-
\begin{equation}\label{80}
f=2^3\cdot7,~2(1+i)17,~2^2\cdot 3\cdot 11,~ 2^3\cdot 17,~ 2(1+i)
\cdot 3\cdot 13,~2^3\cdot 41.
\end{equation}
Tables I and II give the values of $c_p^+(E).p^{-g_{E/\Q}}$ mod $p$
in the two ranges $p <\,1000$ and $29000 \,<p <\, 30000$, with $p$ congruent to 1 mod 4, for our six curves.
We have placed a * in the Tables in the relevant column whenever $p$ divides $D$.

\medskip

Finally, we must  point out that the values of $ c_p^+(E).p^{-g_{E/\Q}}$  \, mod $p$ given in Tables 1
and 2 of \cite{CZS} should be corrected by  multiplying each entry with the factor $4$.  This is because the formulae
(24) and (53) of \cite{CZS} are incorrect, and only become valid if the factor $w$  occurring
in each  formulae is omitted (and $w=4$ when $K=\Q(i)$). Of course, the correct values are given in Tables 1 and 2
below.

\newpage

\begin{center}
\tablefirsthead{%

 \hline

\multicolumn{7}{|c|}{Table I: table  of
$c_p^+(E)/p^{\rm{Rank}(E)}\mod p \text{ for
}p<1000\text{ and }p\equiv 1\mod 4 $.}\\

 \hline$p$& case $82$&case $-33$& case
$-34$&case $-39$&case $17$&case $-14$\\\hline }

\tablehead{\hline \multicolumn{7}{|c|}{Table I: table  of
$c_p^+(E)/p^{\rm{Rank}(E)}\mod p \text{ for
}p<1000\text{ and }p\equiv 1\mod 4 $.}\\

\hline$p$& case $82$&case $-33$& case $-34$&case $-39$&case
$17$&case $-14$\\\hline }

\tabletail{\hline}

 \tablelasttail{\hline}
\begin{supertabular}{|r@{\hspace{15pt}}|r@{\hspace{15pt}}|r@{\hspace{15pt}}|r@{\hspace{15pt}}|r@{\hspace{15pt}}|r@{\hspace{15pt}}|r@{\hspace{15pt}}|}
$ 5 $&$ 3 $&$ 4 $&$ 4 $&$ 1 $&$ 2 $&$ 1 $\\
 $ 13 $&$ 4 $&$ 10 $&$ 9 $&$ * $&$ 6 $&$ 3 $\\
 $ 17 $&$ 14 $&$ 13 $&$ * $&$ 0 $&$ * $&$ 11 $\\
 $ 29 $&$ 20 $&$ 8 $&$ 2 $&$ 20 $&$ 1 $&$ 0 $\\
 $ 37 $&$ 21 $&$ 12 $&$ 26 $&$ 16 $&$ 6 $&$ 36 $\\
 $ 41 $&$ * $&$ 37 $&$ 7 $&$ 37 $&$ 34 $&$ 7 $\\
 $ 53 $&$ 1 $&$ 10 $&$ 14 $&$ 8 $&$ 21 $&$ 9 $\\
 $ 61 $&$ 11 $&$ 1 $&$ 33 $&$ 18 $&$ 43 $&$ 57 $\\
 $ 73 $&$ 42 $&$ 15 $&$ 40 $&$ 44 $&$ 31 $&$ 5 $\\
 $ 89 $&$ 77 $&$ 13 $&$ 76 $&$ 74 $&$ 84 $&$ 82 $\\
 $ 97 $&$ 92 $&$ 57 $&$ 10 $&$ 92 $&$ 41 $&$ 69 $\\
 $ 101 $&$ 31 $&$ 61 $&$ 43 $&$ 92 $&$ 34 $&$ 10 $\\
 $ 109 $&$ 37 $&$ 43 $&$ 31 $&$ 3 $&$ 27 $&$ 54 $\\
 $ 113 $&$ 64 $&$ 29 $&$ 30 $&$ 98 $&$ 31 $&$ 75 $\\
 $ 137 $&$ 6 $&$ 136 $&$ 8 $&$ 7 $&$ 17 $&$ 93 $\\
 $ 149 $&$ 146 $&$ 21 $&$ 89 $&$ 119 $&$ 91 $&$ 146 $\\
 $ 157 $&$ 101 $&$ 100 $&$ 147 $&$ 31 $&$ 109 $&$ 35 $\\
 $ 173 $&$ 158 $&$ 156 $&$ 61 $&$ 3 $&$ 3 $&$ 77 $\\
 $ 181 $&$ 20 $&$ 78 $&$ 45 $&$ 177 $&$ 99 $&$ 85 $\\
 $ 193 $&$ 48 $&$ 72 $&$ 122 $&$ 125 $&$ 74 $&$ 44 $\\
 $ 197 $&$ 110 $&$ 109 $&$ 76 $&$ 14 $&$ 186 $&$ 19 $\\
 $ 229 $&$ 206 $&$ 19 $&$ 51 $&$ 132 $&$ 25 $&$ 207 $\\
 $ 233 $&$ 178 $&$ 140 $&$ 18 $&$ 212 $&$ 136 $&$ 230 $\\
 $ 241 $&$ 114 $&$ 54 $&$ 89 $&$ 108 $&$ 82 $&$ 28 $\\
 $ 257 $&$ 32 $&$ 203 $&$ 41 $&$ 136 $&$ 25 $&$ 36 $\\
 $ 269 $&$ 127 $&$ 261 $&$ 243 $&$ 169 $&$ 214 $&$ 18 $\\
 $ 277 $&$ 50 $&$ 130 $&$ 272 $&$ 267 $&$ 109 $&$ 0 $\\
 $ 281 $&$ 145 $&$ 171 $&$ 50 $&$ 64 $&$ 235 $&$ 147 $\\
 $ 293 $&$ 187 $&$ 232 $&$ 49 $&$ 215 $&$ 121 $&$ 239 $\\
 $ 313 $&$ 65 $&$ 72 $&$ 179 $&$ 123 $&$ 276 $&$ 41 $\\
 $ 317 $&$ 212 $&$ 108 $&$ 14 $&$ 137 $&$ 314 $&$ 130 $\\
 $ 337 $&$ 241 $&$ 104 $&$ 80 $&$ 35 $&$ 76 $&$ 267 $\\
 $ 349 $&$ 114 $&$ 336 $&$ 344 $&$ 85 $&$ 103 $&$ 5 $\\
 $ 353 $&$ 61 $&$ 104 $&$ 78 $&$ 288 $&$ 219 $&$ 60 $\\
 $ 373 $&$ 97 $&$ 195 $&$ 315 $&$ 216 $&$ 300 $&$ 198 $\\
 $ 389 $&$ 381 $&$ 13 $&$ 32 $&$ 273 $&$ 250 $&$ 33 $\\
 $ 397 $&$ 15 $&$ 205 $&$ 279 $&$ 59 $&$ 312 $&$ 272 $\\
 $ 401 $&$ 394 $&$ 197 $&$ 308 $&$ 22 $&$ 193 $&$ 157 $\\
 $ 409 $&$ 255 $&$ 138 $&$ 95 $&$ 70 $&$ 44 $&$ 25 $\\
 $ 421 $&$ 92 $&$ 97 $&$ 369 $&$ 199 $&$ 187 $&$ 134 $\\
 $ 433 $&$ 306 $&$ 155 $&$ 417 $&$ 37 $&$ 429 $&$ 175 $\\
 $ 449 $&$ 11 $&$ 32 $&$ 345 $&$ 337 $&$ 345 $&$ 111 $\\
 $ 457 $&$ 178 $&$ 369 $&$ 159 $&$ 60 $&$ 238 $&$ 133 $\\
 $ 461 $&$ 205 $&$ 304 $&$ 250 $&$ 12 $&$ 71 $&$ 148 $\\
 $ 509 $&$ 335 $&$ 500 $&$ 411 $&$ 208 $&$ 412 $&$ 101 $\\
 $ 521 $&$ 105 $&$ 469 $&$ 172 $&$ 382 $&$ 424 $&$ 129 $\\
 $ 541 $&$ 12 $&$ 18 $&$ 162 $&$ 59 $&$ 132 $&$ 65 $\\
 $ 557 $&$ 229 $&$ 107 $&$ 38 $&$ 216 $&$ 547 $&$ 336 $\\
 $ 569 $&$ 473 $&$ 274 $&$ 566 $&$ 237 $&$ 554 $&$ 267 $\\
 $ 577 $&$ 71 $&$ 506 $&$ 0 $&$ 574 $&$ 156 $&$ 271 $\\
 $ 593 $&$ 505 $&$ 155 $&$ 524 $&$ 454 $&$ 313 $&$ 72 $\\
 $ 601 $&$ 597 $&$ 350 $&$ 515 $&$ 491 $&$ 290 $&$ 229 $\\
 $ 613 $&$ 292 $&$ 311 $&$ 363 $&$ 75 $&$ 490 $&$ 521 $\\
 $ 617 $&$ 408 $&$ 187 $&$ 188 $&$ 206 $&$ 532 $&$ 293 $\\
 $ 641 $&$ 388 $&$ 548 $&$ 186 $&$ 269 $&$ 499 $&$ 33 $\\
 $ 653 $&$ 67 $&$ 343 $&$ 170 $&$ 332 $&$ 384 $&$ 201 $\\
 $ 661 $&$ 642 $&$ 125 $&$ 382 $&$ 176 $&$ 80 $&$ 659 $\\
 $ 673 $&$ 620 $&$ 407 $&$ 102 $&$ 73 $&$ 501 $&$ 115 $\\
 $ 677 $&$ 492 $&$ 583 $&$ 22 $&$ 68 $&$ 651 $&$ 63 $\\
 $ 701 $&$ 682 $&$ 180 $&$ 211 $&$ 527 $&$ 420 $&$ 380 $\\
 $ 709 $&$ 333 $&$ 79 $&$ 707 $&$ 595 $&$ 330 $&$ 432 $\\
 $ 733 $&$ 398 $&$ 544 $&$ 276 $&$ 16 $&$ 307 $&$ 382 $\\
 $ 757 $&$ 577 $&$ 537 $&$ 519 $&$ 150 $&$ 671 $&$ 417 $\\
 $ 761 $&$ 172 $&$ 445 $&$ 554 $&$ 428 $&$ 691 $&$ 101 $\\
 $ 769 $&$ 208 $&$ 448 $&$ 45 $&$ 652 $&$ 697 $&$ 603 $\\
 $ 773 $&$ 264 $&$ 270 $&$ 685 $&$ 439 $&$ 52 $&$ 703 $\\
 $ 797 $&$ 550 $&$ 603 $&$ 33 $&$ 92 $&$ 492 $&$ 372 $\\
 $ 809 $&$ 267 $&$ 606 $&$ 196 $&$ 12 $&$ 154 $&$ 39 $\\
 $ 821 $&$ 638 $&$ 728 $&$ 789 $&$ 45 $&$ 24 $&$ 567 $\\
 $ 829 $&$ 799 $&$ 825 $&$ 430 $&$ 418 $&$ 93 $&$ 805 $\\
 $ 853 $&$ 85 $&$ 476 $&$ 473 $&$ 579 $&$ 192 $&$ 834 $\\
 $ 857 $&$ 13 $&$ 802 $&$ 25 $&$ 5 $&$ 20 $&$ 294 $\\
 $ 877 $&$ 505 $&$ 353 $&$ 734 $&$ 828 $&$ 528 $&$ 658 $\\
 $ 881 $&$ 373 $&$ 427 $&$ 355 $&$ 267 $&$ 328 $&$ 602 $\\
 $ 929 $&$ 775 $&$ 430 $&$ 923 $&$ 9 $&$ 593 $&$ 277 $\\
 $ 937 $&$ 795 $&$ 456 $&$ 829 $&$ 904 $&$ 427 $&$ 400 $\\
 $ 941 $&$ 5 $&$ 594 $&$ 645 $&$ 558 $&$ 71 $&$ 686 $\\
 $ 953 $&$ 35 $&$ 633 $&$ 133 $&$ 819 $&$ 317 $&$ 605 $\\
 $ 977 $&$ 587 $&$ 653 $&$ 428 $&$ 915 $&$ 238 $&$ 392 $\\
 $ 997 $&$ 839 $&$ 563 $&$ 408 $&$ 882 $&$ 211 $&$ 607 $\\

\end{supertabular}

\end{center}

\begin{center}
\tablefirsthead{%

 \hline

\multicolumn{7}{|c|}{Table II: table  of
$c_p^+(E)/p^{\rm{Rank}(E)}\mod p  \text{ for
}29000<p<30000\text{ and }p\equiv 1\mod 4 $.}\\

 \hline$p$& case $82$&case $-33$& case $-34$&case $-39$&case $17$&case $-14$\\\hline }

\tablehead{\hline \multicolumn{7}{|c|}{Table II: table  of
$c_p^+(E)/p^{\rm{Rank}(E)}\mod p  \text{ for
}29000<p<30000\text{ and }p\equiv 1\mod 4 $.}\\

\hline$p$& case $82$&case $-33$& case $-34$&case $-39$&case
$17$&case $-14$\\\hline }

\tabletail{\hline}

 \tablelasttail{\hline}
\begin{supertabular}{|r@{\hspace{15pt}}|r@{\hspace{15pt}}|r@{\hspace{15pt}}|r@{\hspace{15pt}}|r@{\hspace{15pt}}|r@{\hspace{15pt}}|r@{\hspace{15pt}}|}

$ 29009 $&$ 25650 $&$ 23127 $&$ 2319 $&$ 23907 $&$ 2335 $&$ 22753
$\\

 $ 29017 $&$ 2820 $&$ 18581 $&$ 20044 $&$ 17655 $&$ 28801 $&$ 10064 $
\\
 $ 29021 $&$ 11083 $&$ 5690 $&$ 7203 $&$ 16719 $&$ 21473 $&$ 28871
$\\
 $ 29033 $&$ 28457 $&$ 24461 $&$ 8759 $&$ 24782 $&$ 16202 $&$ 22261 $
\\
 $ 29077 $&$ 14804 $&$ 11706 $&$ 10879 $&$ 27218 $&$ 21375 $&$ 26584
$\\
 $ 29101 $&$ 8296 $&$ 17887 $&$ 11405 $&$ 3758 $&$ 4867 $&$ 19197 $\\

 $ 29129 $&$ 25398 $&$ 22104 $&$ 1256 $&$ 15983 $&$ 27700 $&$ 14477 $
\\
 $ 29137 $&$ 11967 $&$ 8087 $&$ 14839 $&$ 2323 $&$ 21504 $&$ 11674
$\\
 $ 29153 $&$ 21169 $&$ 10422 $&$ 27377 $&$ 3436 $&$ 15374 $&$ 11577 $
\\
 $ 29173 $&$ 9596 $&$ 27126 $&$ 10426 $&$ 14543 $&$ 14146 $&$ 15233 $\
\\
 $ 29201 $&$ 27808 $&$ 18751 $&$ 1663 $&$ 26052 $&$ 10329 $&$ 10761 $
\\
 $ 29209 $&$ 17480 $&$ 25574 $&$ 26288 $&$ 23118 $&$ 15890 $&$ 21215
$\\
 $ 29221 $&$ 10441 $&$ 4222 $&$ 3439 $&$ 15921 $&$ 20174 $&$ 2642 $\\

 $ 29269 $&$ 4577 $&$ 8091 $&$ 18622 $&$ 2602 $&$ 8749 $&$ 14748 $\\
 $ 29297 $&$ 302 $&$ 24353 $&$ 11929 $&$ 10928 $&$ 6390 $&$ 24694 $\\

 $ 29333 $&$ 15913 $&$ 6596 $&$ 6496 $&$ 24112 $&$ 25483 $&$ 17739
$\\
 $ 29389 $&$ 28936 $&$ 7811 $&$ 3473 $&$ 19903 $&$ 16483 $&$ 4201 $\\

 $ 29401 $&$ 3675 $&$ 14554 $&$ 27879 $&$ 11800 $&$ 20475 $&$ 13536 $
\\
 $ 29429 $&$ 6780 $&$ 11389 $&$ 26246 $&$ 167 $&$ 24609 $&$ 24452 $\\

 $ 29437 $&$ 17070 $&$ 14093 $&$ 16470 $&$ 28787 $&$ 6001 $&$ 12774 $
\\
 $ 29453 $&$ 57 $&$ 21091 $&$ 12802 $&$ 17808 $&$ 14024 $&$ 15022 $\\

 $ 29473 $&$ 11341 $&$ 8907 $&$ 4052 $&$ 12428 $&$ 27500 $&$ 341 $\\
 $ 29501 $&$ 20384 $&$ 24356 $&$ 12630 $&$ 13100 $&$ 25119 $&$ 13208
$\\
 $ 29537 $&$ 12361 $&$ 4242 $&$ 3350 $&$ 3811 $&$ 3605 $&$ 15245 $\\
 $ 29569 $&$ 2147 $&$ 16852 $&$ 14648 $&$ 3185 $&$ 21444 $&$ 18834
$\\
 $ 29573 $&$ 14394 $&$ 22250 $&$ 1696 $&$ 13205 $&$ 25891 $&$ 6662
$\\
 $ 29581 $&$ 8560 $&$ 16558 $&$ 21191 $&$ 534 $&$ 29153 $&$ 10018 $\\

 $ 29629 $&$ 28337 $&$ 18609 $&$ 8370 $&$ 14834 $&$ 26381 $&$ 5673
$\\
 $ 29633 $&$ 14670 $&$ 5698 $&$ 22556 $&$ 11543 $&$ 19394 $&$ 27146 $
\\
 $ 29641 $&$ 5608 $&$ 20521 $&$ 12418 $&$ 24522 $&$ 7262 $&$ 12119
$\\
 $ 29669 $&$ 20520 $&$ 7822 $&$ 26359 $&$ 26038 $&$ 3437 $&$ 14200
$\\
 $ 29717 $&$ 26628 $&$ 26613 $&$ 16105 $&$ 29029 $&$ 12500 $&$ 20097
$\\
 $ 29741 $&$ 581 $&$ 24054 $&$ 17026 $&$ 5772 $&$ 12560 $&$ 4084 $\\
 $ 29753 $&$ 4788 $&$ 6974 $&$ 14012 $&$ 12099 $&$ 373 $&$ 27858 $\\
 $ 29761 $&$ 10181 $&$ 12106 $&$ 5979 $&$ 5036 $&$ 5256 $&$ 10290 $\\

 $ 29789 $&$ 17699 $&$ 295 $&$ 8422 $&$ 22502 $&$ 17747 $&$ 26607 $\\

 $ 29833 $&$ 14126 $&$ 4510 $&$ 9716 $&$ 481 $&$ 22479 $&$ 4411 $\\
 $ 29837 $&$ 29016 $&$ 20966 $&$ 4124 $&$ 17704 $&$ 5190 $&$ 18103
$\\
 $ 29873 $&$ 24668 $&$ 4490 $&$ 5529 $&$ 5739 $&$ 17532 $&$ 25844 $\\

 $ 29881 $&$ 27804 $&$ 18837 $&$ 27388 $&$ 22372 $&$ 18039 $&$ 23336
$\\
 $ 29917 $&$ 8462 $&$ 17177 $&$ 1116 $&$ 12357 $&$ 3665 $&$ 6666 $\\
 $ 29921 $&$ 14442 $&$ 20982 $&$ 6729 $&$ 14884 $&$ 8950 $&$ 16489
$\\
 $ 29989 $&$ 11184 $&$ 7771 $&$ 6327 $&$ 25027 $&$ 29511 $&$ 22258
$\\

\end{supertabular}
\end{center}

\medskip

Finally, we consider the special primes. The primes $p = 19$ and $p=277$ for the curve with $D=-14$
are already discussed in \cite{CZS}.  For $D = -39$ and  $p=17$, we have
that
$$c_p^+(E)=3^{11}5^{2}7^2 13^4 17^{ 3} 11\cdot 163\cdot 428532544446776087.$$
For $D=-34$,  and  $p=577$, $c_p^+(E)$  lies in the interval $[10^{2289},10^{2290}]$, and is a
very large number, which we have not succeeded in factoring.  However, our calculations show that
 $$c_p^+(E)\equiv 69\cdot 577^3\mod 577^4.$$
Thus, by Theorem 2.2 of \cite{CZS}, $\sha(E/\Q)(p)$ is finite in both cases.

\bigskip

\noindent  J. Coates\hb
\noindent  Emmanuel College\hb
\noindent Cambridge CB2 3AP, England.\hb
\noindent  Department of Mathematics,\hb
POSTECH, Pohang 790-784, Korea.\hb
\noindent e-mail j.h.coates@dpmms.cam.ac

\medskip

\noindent  Z. Liang\hb
\noindent School of Mathematical Sciences,\hb
\noindent Capital Normal University,\hb
\noindent Xisanhuanbeilu 105, Haidan District,\hb
\noindent Beijing, China.\hb
\noindent e-mail liangzhb@gmail.com

\medskip

\noindent R. Sujatha\hb
\noindent School of Mathematics, TIFR,\hb
\noindent Homi Bhabha Road, Colaba,\hb
\noindent Mumbai 40005, India.\hb
\noindent e-mail sujatha@math.tifr.res.in


\begin{thebibliography}{99}
\bibitem{CZS}{\sc J. Coates, Z. Liang, R. Sujatha}, \newblock
{\it The Tate-Shafarevich group for elliptic curves with complex
  multiplication},
Jour. Alg. {\bf 322} (2009), 657--674.
\bibitem{K}{\sc N. Katz}, \newblock
{\it  Divisibilities, congruences, and Cartier duality},
\newblock  Jour. Fac. Sci. Univ. Tokyo {\bf 28} (1982), 667--678.
\bibitem{CW}{\sc J. Coates, A. Wiles},
\newblock{\it On $p$-adic  $L$-functions and elliptic units},
J. Austral. Math. Soc. {\bf 26} (1978), 1-25.
\bibitem{CS}{\sc J. Coates, R. Sujatha},
\newblock{\it Elliptic curves with complex multiplication and
  $L$-values}, book in preparation.
\bibitem{GS}{\sc C. Goldstein, N. Schappacher},
\newblock{\it Series d'Eisenstein e fonctions $L$ des courbes elliptiques
a multiplication complexe}, Crelle J. {\bf 327} (1981), 184-218.
\bibitem{LMS}{\sc J. Coates}, \newblock{\it Elliptic curves with
    complex multiplication and Iwasawa theory},
 Bull. London
Math. Soc. {\bf 23} (1991) 321-350.




\end{thebibliography}
\end{document}